\documentclass[1p,12pt]{elsarticle}

\usepackage{amssymb}
\usepackage{amsmath}
\usepackage{graphicx}
\usepackage{ctable}

\journal{Applied Mathematics and Computation}
\begin{document}

\begin{frontmatter}

\title{Arithmetic properties of homogeneous Hilbert curves}

\author[phys]{E. Estevez-Rams\corref{cor1}}
\ead{estevez@imre.oc.uh.cu}

\author[baf]{I. Brito-Reyes}

\cortext[cor1]{Corresponding author at: Instituto de Ciencias y Tecnolog\'ia de Materiales, University of Havana (IMRE), San Lazaro y L. CP 10400. La Habana. Cuba. t: (+537) 8705707}

\address[phys]{Physics Faculty-Instituto de Ciencias y Tecnolog\'ia de Materiales (IMRE), University of Havana, San Lazaro y L. CP 10400. La Habana. Cuba.}

\address[baf]{Universidad de las Ciencias Inform\'aticas (UCI), Carretera a San Antonio, Boyeros. La Habana. Cuba.}

\begin{abstract}
Properties of the recently reported homogeneous Hilbert curves are deduced and reported. The nature of the affine transformations involved in the construction of the Hilbert curves is explored. The analytical representation of proper and improper Hilbert curves is obtained. The one-to-one mapping between two Hilbert curves is deduced. Recursive relation of Hilbert curves is reported.  
\end{abstract}

\begin{keyword}
Hilbert curve \sep Space filling curves  
\end{keyword}


\date{\today}
\end{frontmatter}

\section{Introduction}

Hilbert was the first to propose an algorithmic procedure for the generation of a space filling curves \cite{sagan94}. Space filling curves were introduced by Peano and was originally defined as a surjective mapping between the unit interval $I$ onto the unit square $Q$. It can be easily generalized to a surjective mapping $\mathbb{R}\longrightarrow\mathbb{R}^{d}$ from the one-dimensional space onto the $d$-dimensional one. Originally considered an oddity and even named ``crinkly'' curves \cite{moore00}, they have found an increasing number of applications in diverse fields \cite{chen05,chen11,bially69,songa02,liang08,anders09}.

Hilbert curves mapping are based on the iterative application of affine transformations to a starting mapping \cite{sagan94}. In the starting mapping the unit interval is partition into four disjoint, equal length, subintervals and put into correspondence with a four disjoint, equal area, subsquares partition of the unit square. At each step, each subinterval and corresponding subsquare, are considered as an original interval and square and the affine transformation is then applied over them. The affine transformation must be so, that continuity is preserved, and two adjacent subintervals are mapped into two adjacent subsquares. The curve is uniquely defined (up to a rotation or reflection) by fixing the mapping of the initial and final subintervals.

Shortly after Hilbert, Moore reported another Hilbert type curve that has subsequently been called by his name \cite{moore00}. Liu \cite{liu04} reported another additional four curves of Hilbert type. Recently, P\'erez-Davidenko et al. \cite{estevez_may13} completed the set of homogeneous Hilbert curves reporting additionally six curves. It was shown that this twelve curves exhaust all possible Hilbert curves in two dimensions made by the application of only one set (homogeneity condition) of affine transformation at each iterative step.   

In this contribution we explore the arithmetic properties of the twelve homogeneous Hilbert curves (HHC) and the relation between them.

\section{Homogeneous Hilbert Curves}\label{sec:HHC}

\begin{figure} \centering\small
\includegraphics[scale=1]{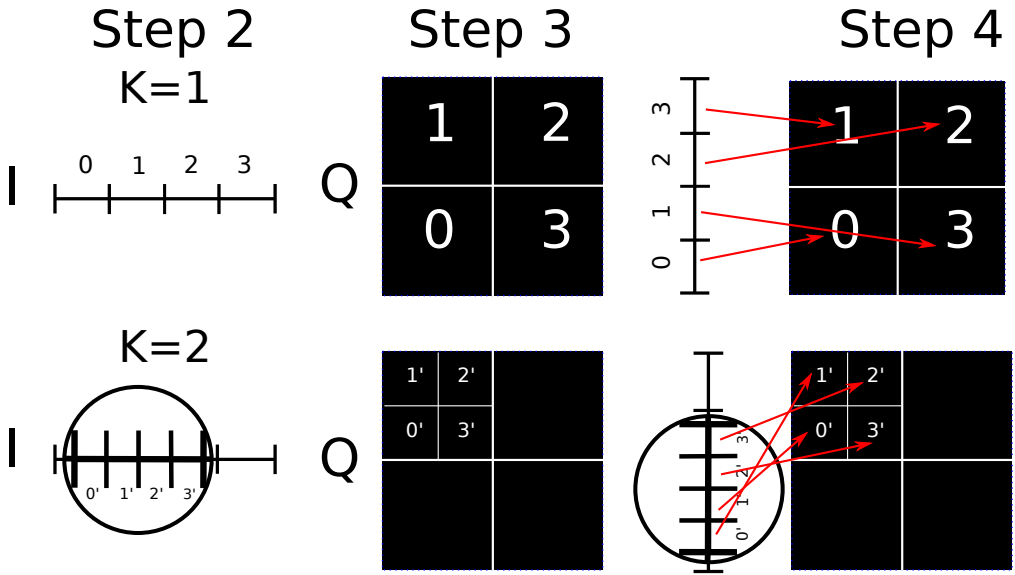}
\caption{Hilbert construction.
}\label{fig:algol}
\end{figure}

The rationale behind Hilbert recursive construction of a space filling curves can be described by the following algorithm (Figure \ref{fig:algol}):

\begin{enumerate}
 \item $I_{0.}=\{r|0 \leq r \leq 1\}$ and $Q_{0.}=\{(x,y)| 0 \leq x \leq 1, 0\leq y \leq 1\}$, $k=0$.
 \item Partition each $I_{0.q_{0}q_{1}\ldots q_{k}}$ into four congruent, disjoint subintervals $I_{0.q_{0}q_{1}\ldots q_{k}}=\{I_{0.q_{0}q_{1}\ldots q_{k}0},I_{0.q_{0}q_{1}\dots q_{k}1},I_{0.q_{0}q_{1}\ldots q_{k}2},I_{0.q_{0}q_{1}\ldots q_{k}3}\}$ (Fig. \ref{fig:algol}a).\label{steppart}
 \item Partition each $Q_{0.m_{0}m_{1}\ldots m_{k}}$ into four congruent, disjoint subsquares $Q_{0.m_{0}m_{1}\ldots m_{k}}=\{Q_{0.m_{0}m_{1}\ldots m_{k}0},Q_{0.m_{0}m_{1}\ldots m_{k}1},Q_{0.m_{0}m_{1}\dots m_{k}2},Q_{0.m_{0}m_{1}\ldots m_{k}3}\}$.
 \item Make the correspondence $I_{0.q_{0}q_{1}\ldots q_{k}q_{k+1}}\longrightarrow Q_{0.m_{0}m_{1}\ldots m_{k}m_{k+1}}$ such that two consecutive subintervals in  $I_{0.q_{0}q_{1}\ldots q_{k}}$ corresponds to two adjacent subsquares in  $Q_{0.m_{0}m_{1}\ldots m_{k}m_{k+1}}$ (adjacency condition).
 \item k=k+1
 \item go to step \ref{steppart}.
\end{enumerate}

The curve obtained at the $k$ iteration is known as the Hilbert curve of order $k$ and will correspond to the partition of the unit interval and the unit square into $4^k$ subintervals and subsquares, respectively. The algorithm described above does not determine uniquely the space filling curve, in order to do so, boundary condition must be given, that is for all $k$, fix the mapping of the initial and final subinterval to the corresponding subsquares in the partition of $Q$.  Hilbert curves are often referred to the infinite limit $k\longrightarrow\infty$, where the mapping becomes surjective.

Hilbert original curve results from mapping at any order $k$, the first subinterval to the lower left subsquare: $I_{0.000\ldots0}\longrightarrow Q_{0.000\ldots0}$; and the last subinterval to the lower right subsquare: $I_{0.333\ldots3}\longrightarrow Q_{0.333\ldots3}$. By assigning different initial and final mapping five additional curves can be constructed, which have been called by Davidenko et al. proper Hilbert curves \cite{estevez_may13}. The boundary conditions for the complete set of proper Hilbert curve can be seen in Table \ref{tbl:boundaryproper}. The proper Hilbert curve of order $k$ will be given the symbol $_{\nu}H_{k}$, with $\nu=0,1,2\ldots 5$, where $_{0}H_{k}$ and $_{1}H_{k}$ are the Hilbert original curve and the Moore curve \cite{moore00}, respectively. For $\nu=2,3,4,5$ the symbol will correspond to the additional proper curves introduced by Liu \cite{liu04}.

If the reversion operation is introduced, six additional curves called improper in \cite{estevez_may13} can be constructed. The reversion operation swaps the entry and exit point of the curve at a quadrant. For all improper curves, the boundary condition for the order 2 curve is the same that the $_{4}H$ (Liu 5) curve. The boundary conditions for the complete set of improper Hilbert curve of order $k>2$ can be seen in Table \ref{tbl:boundaryimproper}.

Proper curves of order $k$ are constructed from a specific mapping, different for each curve type, over the $_{0}H_{k-1}$ curve. Improper curves, on the other hand, are constructed from a specific mapping, different for each curve type, over the $_{5}H_{k-1}$ (Liu 4) curve.

\begin{figure} \centering\small
\includegraphics[scale=0.75]{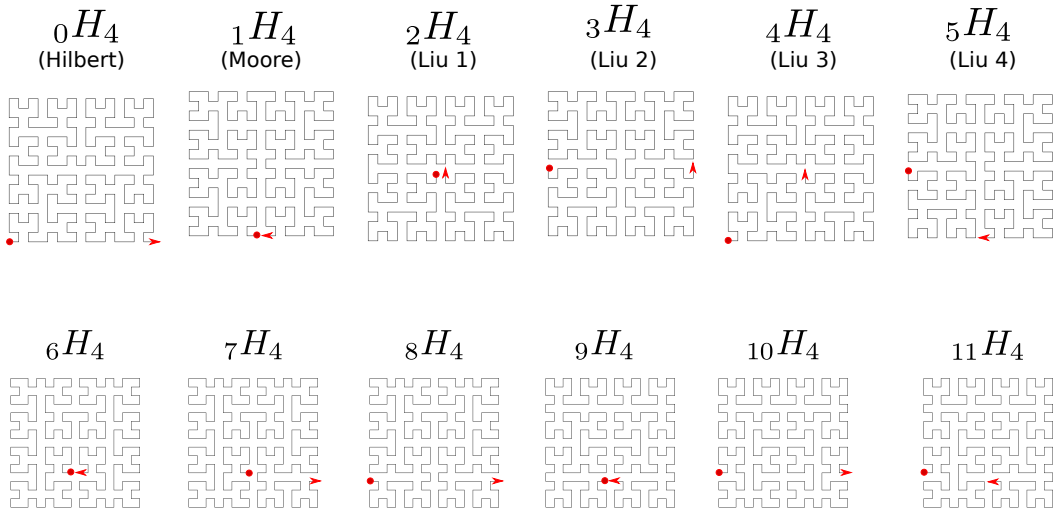}
\caption{The homogeneous Hilbert curves of order 4. The circle signals the entry point, and the arrow the exit point.
}\label{fig:hilbert4}
\end{figure}

Figure \ref{fig:hilbert4} shows the homogeneous Hilbert curves of order $4$.

From a geometric point of view, HHC can be classified by their symmetry: some curves posses  a vertical symmetry line at the middle of the unit square  (m symmetry); they can also be classified according to the nature of the entry and exit points which can lie at a corner, edge or interior subsquare; finally HHC can be closed if the entry and exit point lies at adjacent subsquares (Table \ref{tbl:prop}). 

Homogeneous Hilbert curves can also be described as generated iteratively by the application of certain sets of affine transformations $_{\nu}p_{i}$ that maps one subinterval to the $i$th quadrant. The affine transformations operator has the general form 
\begin{equation}
p\left (\begin{array}{l}x \\ y\end{array}\right )=\frac{1}{2}\mathbf{U}\cdot \left (\begin{array}{l}x \\ y\end{array}\right )+\frac{1}{2}\mathbf{t}=[U,t]_{\frac{1}{2}}\left (\begin{array}{l}x \\ y\end{array}\right ),\label{eq:poperator}
\end{equation}
 where $\mathbf{U}$ is  a rotation, given by a $2\times2$ orthogonal matrix, and $\mathbf{t}$ is a translation factor, given by a vector $(x,y)$. Each $_{\nu}H_{k}$ is constructed by using a set of four $_{\nu}p_{i}\; i=0,1,2,3$ transformation one for each quadrant, that act ($\otimes$) over the coordinates of the a $k-1$-order curve. For the proper curves, the operation can be represented as 

\begin{equation}
_{\nu}H_{n} =_{\nu}H\otimes\;_{0}H_{n-1},\label{eq:pproper}
\end{equation}

while for the improper curves

\begin{equation}\label{eq:pimproper}
_{\nu}H_{n}  =_{\nu}H\otimes\;_{5}H_{n-1}.
\end{equation}

Table \ref{tbl:affine} gives the sets of affine transformations. 

\section{Group properties of the rotation parts of the affine transformations}\label{sec:group}

The rotation matrices of the affine transformations form a group. The multiplication table is shown in table \ref{tbl:multtable}. The group is isomorphic with the planar point group 4mm \cite{ITA} with generators $\{U_{R}, U_{H}\}$. The group structure exhibits three subgroups of order 4: $\{U_{I},-U_{I},U_{R},-U_{R}\}$ and $\{U_{I},-U_{I},U_{H},-U_{H}\}$ isomorphic to 2mm, $\{U_{I},-U_{I},U_{V},-U_{V}\}$ isomorphic to the planar point group 4; five cyclic subgroups of order two: $\{U_{I}, U_{R}\}$, $\{U_{I}, -U_{R}\}$, $\{U_{I}, -U_{I}\}$, $\{U_{I}, U_{H}\}$, $\{U_{I}, -U_{H}\}$.

\section{Arithmetic representation of HHC}\label{sec:arithmetic}

We generalize the procedure for an analytic representation of the original Hilbert curve $_{0}H$ given by Sagan \cite{sagan94} to include all HHC.

Consider the quaternary representation of $t\in [0,1]$:

\begin{equation}
t=\frac{q_{1}}{4}+\frac{q_{2}}{4^{2}}+\frac{q_{3}}{4^{3}}+\ldots+\frac{q_{k}}{4^{k}}\equiv0.q_{1}q_{2}q_{3}\ldots q_{k}
\end{equation}
 
\begin{figure} \centering\small
\includegraphics[scale=0.75]{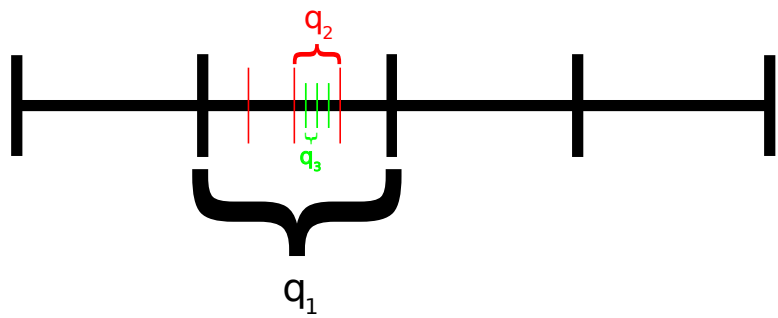}
\caption{The quaternary decomposition of a number.
}\label{fig:quater}
\end{figure}

$q_{s}=0,1,2,3$. By construction $t$ lies in the $q_{1}+1$-th subinterval of the first partition of the unit interval; the $q_{2}+1$-th subinterval of the second partition; the $q_{3}+1$-th subinterval of the third partition and so on, up to the $q_{k}+1$-th subinterval of the kth-partition (Figure \ref{fig:quater}).

\subsection{Proper HHC}\label{sec:properhhc}

For the proper Hilbert curves, according to equation (\ref{eq:pproper}), under the $_{\nu}H_{k}$ ($\nu=0,1,\ldots 5$) mapping $f^{(k)}_{\nu}: I\rightarrow Q$, 

\begin{equation}
 f^{(k)}_{\nu}(t)=\,_{\nu}p_{q_{1}} \,_{0}p_{q_{2}} \,_{0}p_{q_{3}}\ldots \,_{0}p_{q_{k}}\Omega\label{eq:pop}
\end{equation}

where $\Omega$, is the starting configuration, we can take $\Omega=\{(1/2,1/2)\}$. From equation (\ref{eq:poperator})
 
\begin{equation}
\begin{array}{rl}
 f^{(k)}_{\nu}(t)=&\,_{\nu}p_{q_{1}} \,_{0}p_{q_{2}} \,_{0}p_{q_{3}}\ldots \,_{0}p_{q_{k-1}} [ \frac{1}{2} \,_{0}U_{q_{k}}\Omega+\frac{1}{2}\,_{0}t_{q_{k}}  ]=\\
&=\,_{\nu}p_{q_{1}} \,_{0}p_{q_{2}} \,_{0}p_{q_{3}}\ldots \,_{0}p_{q_{k-2}} [\frac{1}{2}\,_{0}U_{q_{k-1}}[ \frac{1}{2} \,_{0}U_{q_{k}}\Omega+\frac{1}{2}\,_{0}t_{q_{k}}  ]+\frac{1}{2}\,_{0}t_{q_{k-1}}]=\\
&=\,_{\nu}p_{q_{1}} \,_{0}p_{q_{2}} \,_{0}p_{q_{3}}\ldots \,_{0}p_{q_{k-2}} [\frac{1}{2^{2}}\,_{0}U_{q_{k-1}}\,_{0}U_{q_{k}}\Omega+\frac{1}{2^{2}}\,_{0}U_{q_{k-1}}\,_{0}t_{q_{k}}+\frac{1}{2}\,_{0}t_{q_{k-1}}]=\\
&=\displaystyle \left ( \frac{1}{2^{k}}\right )\,_{\nu}U_{q_{1}}\,_{0}U_{q_{2}}\ldots \,_{0}U_{q_{k}}\Omega+\\
&+\sum_{j=2}^{k}\left ( \frac{1}{2^{j}}\right )\,_{\nu}U_{q_{1}}\,_{0}U_{q_{2}}\ldots _{0}U_{q_{k-1}}\,_{0}t_{q_{j}}+\frac{1}{2}\,_{\nu}t_{q_{1}}.\label{eq:ftk}
\end{array}
\end{equation}

Equation (\ref{eq:ftk}) gives the arithmetic representation of a Hilbert curve of order $k$: $_{\nu}H_{k}$.

In the limit $k\longrightarrow \infty$, the first term in equation (\ref{eq:ftk}) tends to zero and we arrive to the following arithmetic representation of the space filling $_{\nu}H$ curve:
 \begin{equation}
 f_{\nu}(t)=\displaystyle \sum_{j=2}^{k}\left ( \frac{1}{2^{j}}\right )\,_{\nu}U_{q_{1}}\,_{0}U_{q_{2}}\ldots _{0}U_{q_{k-1}}\,_{0}t_{q_{j}}+\frac{1}{2}\,_{\nu}t_{q_{1}}.\label{eq:ft}
\end{equation}
 
From the multiplication table (Table \ref{tbl:multtable}), $-U_{I}=-U_{R}\cdot U_{R}$, where $U_{I}$ is the identity matrix, which simplifies equation (\ref{eq:ftk}) to
\begin{eqnarray*}
 f^{(k)}_{\nu}(t)=&\displaystyle \left ( \frac{1}{2^{k}}\right ) \left ( -1\right )^{\#_{3}(2, k)}\,_{\nu}U_{q_{1}}\Omega+\\
&+\sum_{j=2}^{k}\left ( \frac{1}{2^{j}}\right ) \left ( -1 \right )^{\#_{3}(2,j-1)}\,_{\nu}U_{q_{1}} U_{R}^{\#_{03}(2,j-1)}\,_{0}t_{q_{j}}+\frac{1}{2}\,_{\nu}t_{q_{1}},
\end{eqnarray*}
where use have been made of $U_{R}\Omega=U_{R}\left (\begin{array}{c}1/2 \\1/2\end{array}\right )=\Omega$. $\#_{3}(s,t)$ is the number of 3, from $q_{s}$ to $q_{t}$, in the quaternary expansion of $t$ ; correspondingly, $\#_{03}(s,t)$ is the number of 0 and 3 , from $q_{s}$ to $q_{t}$, in the quaternary expansion of $t$:
\begin{equation*}
 \begin{array}{l}
 \#_{3}(s,t)=\displaystyle \frac{1}{6}\sum_{m=s}^{t} q_{m}(q_{m}-2)(q_{m}-1)\\
 \#_{03}(s,t)=\displaystyle t+\frac{1}{2}\sum_{m=s}^{t} q_{m}(q_{m}-3)
\end{array}
\end{equation*}

If $\#_{03}(s,t)$ is even, then $U_{R}^{\#_{03}(s,t)}=U_{I}$ the identity operator, otherwise, $U_{R}^{\#_{03}(s,t)}=U_{R}$.

According to equation (\ref{eq:pop}),
\begin{equation*}
 f^{(k)}_{0}(t)=\,_{0}p_{q_{1}} \,_{0}p_{q_{2}} \,_{0}p_{q_{3}}\ldots \,_{0}p_{q_{k}},\Omega
\end{equation*}
from where
\begin{equation}
\,_{0}p_{q_{1}}^{-1}  f^{(k)}_{0}(t)=\,_{0}p_{q_{2}} \,_{0}p_{q_{3}}\ldots \,_{0}p_{q_{k}}\Omega,\label{eq:pop0}
\end{equation}
comparing equation (\ref{eq:pop0}) with equation (\ref{eq:pop}), results in
\begin{equation*}
 f^{(k)}_{\nu}(t)=\,_{\nu}p_{q_{1}} \,_{0}p_{q_{1}}^{-1}  f^{(k)}_{0}(t).
\end{equation*}
From this last equation, for two Hilbert curves $_{\nu}H_{k}$  and $_{\nu'}H_{k}$
\begin{equation*}
\begin{array}{l}
 f^{(k)}_{0}(t)=\,_{0}p_{q_{1}} \,_{\nu}p_{q_{1}}^{-1}  f^{(k)}_{\nu}(t)\\
 f^{(k)}_{0}(t)=\,_{0}p_{q_{1}} \,_{\nu'}p_{q_{1}}^{-1}  f^{(k)}_{\nu'}(t),
\end{array}
\end{equation*}

and finally we get

\begin{equation}
 f^{(k)}_{\nu}(t)=\,_{\nu}p_{q_{1}} \,_{\nu'}p_{q_{1}}^{-1}  f^{(k)}_{\nu'}(t)\label{eq:popnu}
\end{equation}
which maps one-to-one a point in the proper Hilbert curve $_{\nu}H_{k}$ with its corresponding (in the sense of belonging to the same value $t$) point in another proper Hilbert curve $_{\nu'}H_{k}$.

From equation (\ref{eq:popnu}), two points lying in the same quadrant preserve their distance when changing from one Hilbert curve to another. Indeed, if two points $t$ and $t_{1}$ are mapped to the same quadrant, then in their quaternary decomposition, the share the same $q_{1}$ value, and
\begin{equation*}
\begin{array}{rl}
 \left |f^{(k)}_{\nu}(t)-f^{(k)}_{\nu}(t_{1})\right |=&\left |\,_{\nu}p_{q_{1}} \,_{\nu'}p_{q_{1}}^{-1}  f^{(k)}_{\nu'}(t)-\,_{\nu}p_{q_{1}} \,_{\nu'}p_{q_{1}}^{-1}  f^{(k)}_{\nu'}(t_{1})\right |\\\\
=& \left |\,_{\nu}U_{q_{1}} \,_{\nu'}U_{q_{1}}^{-1} \left [ f^{(k)}_{\nu'}(t)-f^{(k)}_{\nu'}(t_{1}) \right ]\right |\\\\
=& \left |\,_{\nu}U_{q_{1}}\right | \left | \,_{\nu'}U_{q_{1}}^{-1}\right| \left |  f^{(k)}_{\nu'}(t)-f^{(k)}_{\nu'}(t_{1}) \right |\\\\
=&\left |  f^{(k)}_{\nu'}(t)-f^{(k)}_{\nu'}(t_{1}) \right |
\end{array}
\end{equation*}
where in the last step, use have been of the orthogonality of the $_{\nu}U$ in the affine transformations (Table \ref{tbl:affine}).

\subsection{Improper HHC}\label{sec:improperhhc}

According to equation (\ref{eq:pimproper}),
\begin{equation}
\begin{array}{rl}
 f^{(k)}_{\nu}(t)=&\,_{\nu}p_{q_{1}} \,_{5}p_{q_{2}} \,_{0}p_{q_{3}}\ldots \,_{o}p_{q_{k}} \Omega\\\\
=&\,_{\nu}p_{q_{1}}f^{(k-1)}_{5}(0.q_{2}q_{3}\ldots q_{k})\label{eq:opimp}
\end{array}
\end{equation}

which, using equation (\ref{eq:ftk}), leads to  
\begin{equation}
\begin{array}{rl}
 f^{(k)}_{\nu}(t)=&\displaystyle \left ( \frac{1}{2^{k}}\right )\,_{\nu}U_{q_{1}}\,_{5}U_{q_{2}}\,_{0}U_{q_{3}}\ldots \,_{0}U_{q_{k}}\Omega+\\\\
&+\sum_{j=3}^{k}\left ( \frac{1}{2^{j}}\right )\,_{\nu}U_{q_{1}}\,_{5}U_{q_{2}}\,_{0}U_{q_{3}}\ldots _{0}U_{q_{k-1}}\,_{0}t_{q_{j}}+\frac{1}{4}\,_{\nu}U_{q_{1}}\,_{5}t_{q_{2}}+\frac{1}{2}\,_{\nu}t_{q_{1}}.\label{eq:ftkimp}
\end{array}
\end{equation}

Besides equation (\ref{eq:ftkimp}), reversion operation must be taken care of. Depending on which quadrant, the number $t$ is mapped (the integer part of $4t$:  $\lfloor 4 t \rfloor$), the transformation $t\longrightarrow t'$ must be made according to Table \ref{tbl:reverse}.

\section{Recursive relations}

Let the quaternary decomposition of $t\in [0,1]$ be 
\begin{equation*}
t=0.q_{1}q_{2}q_{3}\ldots q_{k},
\end{equation*}
then, it is easily realized that
\begin{equation}
\displaystyle \begin{array}{rcl}&m&\\\frac{t}{4^m}=0.&\overbrace{00\ldots0}&q_{1}q_{2}q_{3}\ldots q_{k}\end{array}.\label{eq:qt4}
\end{equation}

\subsection{Proper HHC}\label{sec:properhhc1}

According to equations (\ref{eq:pop}) and (\ref{eq:qt4}),
\begin{equation}
 f^{(k+1)}_{\nu}(t/4)=\,_{\nu}p_{0}\cdot f^{(k)}_{0}(t),\label{eq:ft1}
\end{equation}
and
\begin{equation}
 f^{(k+m)}_{\nu}(t/4^{m})=\,_{\nu}p_{0}\cdot \,_{0}p_{0}^{m-1}\cdot f^{(k)}_{0}(t).\label{eq:ftm}
\end{equation}
If $m$ is odd, then
\begin{equation}
 \,_{0}p_{0}^{m-1}=\left [ \frac{1}{2^{m-1}}U_{I}, t_{0}\right],
\end{equation}
which results in
 \begin{equation}
\,_{\nu}p_{0}\cdot \,_{0}p_{0}^{m-1}=\left [ \frac{1}{2^{m-1}}U_{I}, \frac{2^{m-1}-1}{2^{m}}\,_{\nu}t_{0}\right]\cdot\,_{\nu}p_{0}.\label{eq:nup0}
\end{equation}
Substituting equation (\ref{eq:nup0}) in equation (\ref{eq:ftm}) and making use of equation (\ref{eq:ft1}) we finally get,
\begin{equation}
\begin{array}{rl}
 f^{(k+2 r+1)}_{\nu}(t/4^{2r+1})=&\left [ \frac{1}{4^{r}}U_{I}, \frac{4^{r}-1}{2^{2r+1}}\,_{\nu}t_{0}\right]\cdot f^{(k+1)}_{0}(t/4)\\\\
=&\frac{1}{4^{r}}f^{(k+1)}_{0}(t/4)+\frac{4^{r}-1}{2^{2r+1}}\,_{\nu}t_{0}
\end{array}\label{eq:ftrecodd}
\end{equation}

If $m$ is even, then
\begin{equation}
 \,_{0}p_{0}^{m-1}=\left [ \frac{1}{2^{m-1}}U_{R}, t_{0}\right],
\end{equation}
which, after some algebraic manipulation and using the multiplication table (Table \ref{tbl:multtable}),  results in
\begin{equation}
\begin{array}{rl}
 f^{(k+2 r)}_{\nu}(t/4^{2r})=&\left [ \frac{(-1)^{\nu}}{2^{2r-1}}U_{R},\frac{1}{2}\,\nu t_{0}-\frac{(-1)^{\nu}}{4^{r}}U_{R}\cdot \,\nu t_{0}\right]\cdot f^{(k+1)}_{0}(t/4)\\\\
=&\frac{(-1)^{\nu}}{2^{2r-1}}U_{R}\cdot f^{(k+1)}_{0}(t/4)+\frac{1}{2}\,\nu t_{0}-\frac{(-1)^{\nu}}{4^{r}}U_{R}\cdot \,\nu t_{0}
\end{array}\label{eq:ftreceven}
\end{equation}

\subsection{Improper HHC}\label{sec:improperhhc1}

The recursive relation for the improper HHC can readily obtained from equation (\ref{eq:opimp}),
\begin{equation}
\begin{array}{rl}
 f_{\nu}^{(k+m)}(t/4^{m})=&\,_{\nu}p_{q_{1}}f_{5}^{(k+m-1)}(t/4^{m-1})\nonumber\\\\
=&\frac{1}{2}\,_{\nu}U_{q_{1}}f_{5}^{(k+m-1)}(t/4^{m-1})+\frac{1}{2}\,_{\nu}t_{q_{1}}
\end{array}
\end{equation}
and using the recursive relation for $f_{5}^{(k+p)}(t/4^{p})$ from equations (\ref{eq:ftrecodd}) or (\ref{eq:ftreceven}).

\section{Concluding Remarks}

Literature has almost exclusively dealt with Hilbert original construction. The recent report of the complete set of homogeneous Hilbert curves widens the set of possible Hilbert type space filling curves for a number of applications which could benefit from the availability of different choices. To set a sound base for the use of HHC, they need to be well understood geometrically and analytically. It is also important to know the relation between HHC and Hilbert original curve.

Open questions still remains, for example, numerical calculation up to high orders seems to suggest that all HHC have the same dilation factor in the infinite iterative limit, which for the original Hilbert curve was proven to be 6 \cite{bauman06}, yet a proof of such conjecture still has not been reported. The tools developed here, could prove useful in such a proof. 

\section{Acknowledgments}

Universidad de la Habana and Universidad de la Ciencias Inform\'aticas are acknowledge for financial support and computational infrastructure.

\bibliographystyle{elsarticle-num}
\bibliography{arithmeticproperties}

\pagebreak

\begin{table}
\caption{The boundary conditions for the complete set of Hilbert proper curves. The fourth column are the boundary conditions for the vertical mirror reflected curve equivalent to the original one.}\label{tbl:boundaryproper} \centering\small
\begin{tabular}{l>{\bfseries}lll}
$\nu$ & Curve      & Boundary conditions & Mirror reflected\\ 
\toprule
0 & Hilbert    & $I_{0.000\ldots0}\longrightarrow Q_{0.000\ldots0}$ &\\
  &            & $I_{0.333\ldots3}\longrightarrow Q_{0.333\ldots3}$ &\\ 
\hline
1 & Moore      & $I_{0.000\ldots0}\longrightarrow Q_{0.033\ldots3}$ &\\
  &            & $I_{0.333\ldots3}\longrightarrow Q_{0.300\ldots0}$ &\\
\hline 
2 & Liu 1      & $I_{0.000\ldots0}\longrightarrow Q_{0.022\ldots2}$ &\\
  &            & $I_{0.333\ldots3}\longrightarrow Q_{0.311\ldots1}$ &\\
\hline 
3 & Liu 2      & $I_{0.000\ldots0}\longrightarrow Q_{0.011\ldots1}$ &\\
  &            & $I_{0.333\ldots3}\longrightarrow Q_{0.322\ldots2}$ &\\ 
\hline
4 & Liu 3      & $I_{0.000\ldots0}\longrightarrow Q_{0.000\ldots0}$ & $I_{0.000\ldots0}\longrightarrow Q_{0.022\ldots2}$\\
  &            & $I_{0.333\ldots3}\longrightarrow Q_{0.311\ldots1}$ & $I_{0.333\ldots3}\longrightarrow Q_{0.333\ldots3}$\\ 
\hline
5 & Liu 4      & $I_{0.000\ldots0}\longrightarrow Q_{0.011\ldots1}$ & $I_{0.000\ldots0}\longrightarrow Q_{0.033\ldots3}$\\
  &            & $I_{0.333\ldots3}\longrightarrow Q_{0.300\ldots0}$ & $I_{0.333\ldots3}\longrightarrow Q_{0.322\ldots2}$\\ 
\bottomrule
\end{tabular}
\end{table}

\begin{table}
\caption{The boundary conditions for the complete set of Hilbert improper curves with order larger than 3. The fourth column are the boundary conditions for the vertical mirror reflected curve equivalent to the original one.}\label{tbl:boundaryimproper} \centering\small
\begin{tabular}{l>{\bfseries}lll}
$\nu$ & Curve      & Boundary conditions & Mirror reflected\\ 
\toprule
6 & I1      & $I_{0.000\ldots0}\longrightarrow Q_{0.023\ldots3}$ &\\
  &         & $I_{0.333\ldots3}\longrightarrow Q_{0.310\ldots0}$ &\\ 
\hline
7 & I2      & $I_{0.000\ldots0}\longrightarrow Q_{0.023\ldots3}$ & $I_{0.000\ldots0}\longrightarrow Q_{0.001\ldots1}$ \\
  &         & $I_{0.333\ldots3}\longrightarrow Q_{0.332\ldots2}$ & $I_{0.333\ldots3}\longrightarrow Q_{0.310\ldots0}$ \\
\hline 
8 & I3      & $I_{0.000\ldots0}\longrightarrow Q_{0.001\ldots1}$ &\\
  &         & $I_{0.333\ldots3}\longrightarrow Q_{0.332\ldots2}$ &\\
\hline 
9 & I4      & $I_{0.000\ldots0}\longrightarrow Q_{0.032\ldots2}$ &\\
  &         & $I_{0.333\ldots3}\longrightarrow Q_{0.301\ldots1}$ &\\ 
\hline
10 & I5     & $I_{0.000\ldots0}\longrightarrow Q_{0.010\ldots0}$ \\
  &         & $I_{0.333\ldots3}\longrightarrow Q_{0.323\ldots3}$ \\ 
\hline
11 & I6     & $I_{0.000\ldots0}\longrightarrow Q_{0.010\ldots0}$ & $I_{0.000\ldots0}\longrightarrow Q_{0.032\ldots2}$\\
  &         & $I_{0.333\ldots3}\longrightarrow Q_{0.301\ldots1}$ & $I_{0.333\ldots3}\longrightarrow Q_{0.323\ldots3}$\\ 
\bottomrule
\end{tabular}
\end{table}

\begin{table}
\caption{Geometric properties of HHC.}\label{tbl:prop}
\begin{tabular}{llll}
\toprule
Curve               & Symm. &  Entry-Exit points         & Closed \\   
\midrule
$_{0}H$ (Hilbert)   & m     & corner-corner  shared edge &        \\
$_{1}H$ (Moore)     & m     & edge-edge shared           & X      \\
$_{2}H$ (Liu1)      & m     & interior-interior          & X      \\
$_{3}H$ (Liu2)      & m     & edge-edge opposed          &        \\
$_{4}H$ (Liu3)      & 1     & corner-interior            &        \\
$_{5}H$ (Liu4)      & 1     & edge-edge  adjacent        &        \\
$_{6}H$ (I1)        & m     & interior-interior          & X      \\
$_{7}H$ (I2)        & 1     & interior-edge              &        \\
$_{8}H$ (I3)        & m     & edge-edge opposed          &        \\
$_{9}H$ (I4)        & m     & interior-interior          & X      \\
$_{10}H$ (I5)       & m     & edge-edge opposed          &        \\
$_{11}H$ (I6)       & 1     & edge-interior              &        \\
\bottomrule 
\end{tabular}
\end{table}

\begin{table}
\caption{Affine transformations $p=[U,\mathbf{t}]_{\frac{1}{2}}$ for the HHC. (The $\frac{1}{2}$ subscript is dropped for succinctness. The overbar means reversion operation, see \cite{estevez_may13} for details).}\centering\small
\begin{tabular}{lllll}
\toprule
Curve      & $q_{0}$ & $q_{1}$ & $q_{2}$ & $q_{3}$ \\
\midrule
$_{0}H$ (Hilbert)    & $\left [ U_{R}, t_{0}\right ]$ & $\left [ U_{I}, t_{1}\right ]$ & $\left [ U_{I}, t_{3}\right ]$ & $\left [ -U_{R}, t_{4}\right ]$ \\ 
$_{1}H$ (Moore)    & $\left [ U_{V}, t_{2}\right ]$ & $\left [ U_{V}, t_{3}\right ]$ & $\left [ -U_{V}, t_{5}\right ]$ & $\left [ -U_{V}, t_{3}\right ]$ \\
$_{2}H$ (Liu 1)    & $\left [ -U_{I}, t_{3}\right ]$ & $\left [ U_{I}, t_{1}\right ]$ & $\left [ U_{I}, t_{3}\right ]$ & $\left [ -U_{I}, t_{4}\right ]$ \\
$_{3}H$ (Liu 2)    & $\left [ U_{H}, t_{1}\right ]$ & $\left [ U_{V}, t_{3}\right ]$ & $\left [ -U_{V}, t_{5}\right ]$ & $\left [U_{H}, t_{3}\right ]$ \\
$_{4}H$ (Liu 3)    & $\left [ U_{R}, t_{0}\right ]$ & $\left [ U_{I}, t_{1}\right ]$ & $\left [ U_{I}, t_{3}\right ]$ & $\left [-U_{I}, t_{4}\right ]$ \\
$_{5}H$ (Liu 4)    & $\left [ U_{H}, t_{1}\right ]$ & $\left [ U_{V}, t_{3}\right ]$ & $\left [ -U_{V}, t_{5}\right ]$ & $\left [-U_{V}, t_{3}\right ]$ \\
$_{6}H$ (I1)    & $\left [ -U_{I}, t_{3}\right ]$ & $\overline{\left [ -U_{H}, t_{3}\right ]}$ & $\left [ U_{I}, t_{3}\right ]$ & $\overline{\left [U_{H}, t_{3}\right ]}$ \\
$_{7}H$ (I2)    & $\left [ -U_{I}, t_{3}\right ]$ & $\overline{\left [ -U_{H}, t_{3}\right ]}$ & $\left [ U_{I}, t_{3}\right ]$ & $\left [-U_{R}, t_{4}\right ]$ \\
$_{8}H$ (I3)    & $\overline{\left [ -U_{V}, t_{1}\right ]}$ & $\overline{\left [ -U_{H}, t_{3}\right ]}$ & $\left [ U_{I}, t_{3}\right ]$ & $\left [-U_{R}, t_{4}\right ]$ \\
$_{9}H$ (I4)    & $\overline{\left [ -U_{R}, t_{3}\right ]}$ & $\left [ U_{V}, t_{3}\right ]$ & $\overline{\left [ U_{R}, t_{3}\right ]}$ & $\left [-U_{V}, t_{3}\right ]$ \\
$_{10}H$ (I5)    & $\left [ U_{H}, t_{1}\right ]$ & $\left [ U_{V}, t_{3}\right ]$ & $\overline{\left [ U_{R}, t_{3}\right ]}$ & $\overline{\left [-U_{I}, t_{4}\right ]}$ \\
$_{11}H$ (I6)    & $\left [ U_{H}, t_{1}\right ]$ & $\left [ U_{V}, t_{3}\right ]$ & $\overline{\left [ U_{R}, t_{3}\right ]}$ & $\left [-U_{V}, t_{3}\right ]$ \\
\bottomrule
\end{tabular}\label{tbl:affine}

\begin{tabular}{l}
\\
\end{tabular}

\begin{tabular}{ll}
$U_{I}=\left ( \begin{array}{cc} 1 & 0  \\ 0 & 1 \end{array} \right )$ &
$U_{R}=\left ( \begin{array}{cc} 0 & 1  \\ 1 & 0 \end{array} \right ) $ \\
$U_{V}=\left ( \begin{array}{cc} 0 & -1  \\ 1 & 0 \end{array} \right )$ &
$U_{H}=\left ( \begin{array}{cc} 1 & 0  \\ 0 & -1 \end{array} \right )$
\end{tabular}

\begin{tabular}{l}
\\
\end{tabular}

\begin{tabular}{lll}
$t_{0}=\left ( \begin{array}{c} 0 \\ 0 \end{array} \right ) $ &
$t_{1}=\left ( \begin{array}{c} 0 \\ 1 \end{array} \right ) $ &
$t_{2}=\left ( \begin{array}{c} 1 \\ 0 \end{array} \right ) $ \\
$t_{3}=\left ( \begin{array}{c} 1 \\ 1 \end{array} \right ) $ &
$t_{4}=\left ( \begin{array}{c} 2 \\ 1 \end{array} \right ) $ &
$t_{5}=\left ( \begin{array}{c} 1 \\ 2 \end{array} \right ) $ 
\end{tabular}
\end{table}

\begin{table}
\caption{The multiplication table for the rotation parts of the affine transformations. The rotation operations are defined in Table \ref{tbl:affine}.}\label{tbl:multtable}\centering\small
\begin{tabular}{r|rrrrrrrr}
        & $U_{I}$  & $U_{R}$  & $-U_{I}$ & $-U_{R}$ & $U_{V}$  & $U_{H}$  & $-U_{V}$ & $-U_{H}$ \\
\hline
$U_{I}$ & $\mathbf{U_{I}}$  & $U_{R}$  & $-U_{I}$ & $-U_{R}$ & $U_{V}$  & $U_{H}$  & $-U_{V}$ & $-U_{H}$ \\
$U_{R}$ & $\mathbf{U_{R}}$  & $U_{I}$  & $-U_{R}$ & $-U_{I}$ & $U_{H}$  & $U_{V}$  & $-U_{H}$ & $-U_{V}$ \\
$-U_{I}$ & $\mathbf{-U_{I}}$ & $-U_{R}$ & $U_{I}$  & $U_{R}$  & $-U_{V}$ & $-U_{H}$ & $U_{V}$  & $U_{H}$ \\
$-U_{R}$ & $\mathbf{-U_{R}}$ & $-U_{I}$ & $U_{R}$  & $U_{I}$  & $-U_{H}$ & $-U_{V}$ & $U_{H}$  & $U_{V}$ \\
$U_{V}$ & $U_{V}$  & $-U_{H}$ & $\mathbf{-U_{V}}$ & $U_{H}$  & $-U_{I}$ & $U_{R}$  & $U_{I}$  & $-U_{R}$ \\
$U_{H}$ & $\mathbf{U_{H}}$  & $-U_{V}$ & $-U_{H}$ & $U_{V}$  & $-U_{R}$ & $U_{I}$  & $U_{R}$  & $-U_{I}$ \\
$-U_{V}$ & $-U_{V}$ & $U_{H}$  & $\mathbf{U_{V}}$  & $-U_{H}$ & $U_{I}$  & $-U_{R}$ & $-U_{I}$ & $U_{R}$ \\
$-U_{H}$ & $\mathbf{-U_{H}}$ & $U_{V}$  & $U_{H}$  & $-U_{V}$  & $U_{R}$  & $-U_{I}$ & $-U_{R}$ & $U_{I}$ 
\end{tabular}
\end{table}

\begin{table}
\caption{Transformation of value $t$ ( $t\longrightarrow t'$) for the improper $_{\nu}H_{k}$, before using equation (\ref{eq:ftkimp}).}\centering\small
\begin{tabular}{l|llll}
\toprule
$\lfloor 4t \rfloor$     & 0 & 1 & 2 & 3 \\
\midrule
$_{6}H$ (I1)    & $t$            & $3/4-1/4^k-t$ & $t$               & $7/4-1/4^k-t$  \\
$_{7}H$ (I2)    & $t$            & $3/4-1/4^k-t$ & $t$               & $t$            \\
$_{8}H$ (I3)    & $1/4-1/4^k-t$  & $3/4-1/4^k-t$ & $t$               & $t$            \\
$_{9}H$ (I4)    & $1/4-1/4^k-t$  & $t$           & $5/4-1/4^k-t$     & $t$            \\
$_{10}H$ (I5)   & $t$            & $t$           & $5/4-1/4^k-t$     & $7/4-1/4^k-t$  \\
$_{11}H$ (I6)   & $t$            & $t$           & $5/4-1/4^k-t$     & $t$            \\
\bottomrule
\end{tabular}\label{tbl:reverse}
\end{table}

\end{document}